\documentclass[10pt,a4]{amsart}
\usepackage{verbatim}

\usepackage{latexsym}

\newif\ifpdf
\ifx\pdfoutput\undefined
\pdffalse 
\else
\pdfoutput=1 
\pdftrue
\fi

\ifpdf
\usepackage[pdftex]{graphicx}
\else
\usepackage{graphicx}
\fi

\newtheorem{thm}{{\sc Theorem}}[section]

\newtheorem{lem}[thm]{{\sc Lemma}}

\numberwithin{equation}{section}

\newcommand\pa{\partial}
\newcommand\Lap{\Delta}
\newcommand\RR{\mathbb{R}}
\newcommand\cl{\operatorname{cl}}

\newcommand\ubar{\overline{u}}
\newcommand\ep{\epsilon}
\newcommand\phiep{\phi_{\epsilon}}
\newcommand\uep{u_{\epsilon}}
\newcommand\restriction{\, | \, }

\newcommand\be{\begin{equation}}
\newcommand\Imag{\operatorname{Im}}

\title[Upper and lower bounds for Dirichlet eigenfunctions]
{Upper and lower bounds for normal derivatives of Dirichlet eigenfunctions}

\author{Andrew Hassell}
\address{Department of Mathematics, ANU, Canberra, ACT 0200 \\ AUSTRALIA}
\email{hassell@maths.anu.edu.au}

\author{Terence Tao}
\address{Department of Mathematics, UCLA, Los Angeles CA 90095-1555 \\ USA}
\email{tao@math.ucla.edu}

\thanks{A. H.  is supported by an Australian Research Council Fellowship. T. T.  is a Clay Prize Fellow and is supported by the Packard Foundation. }

\begin{document}

\ifpdf
\DeclareGraphicsExtensions{.pdf, .jpg}
\else
\DeclareGraphicsExtensions{.eps, .jpg}
\fi

\begin{abstract}Suppose that $M$ is a compact Riemannian manifold with boundary and $u$ is an $L^2$-normalized Dirichlet eigenfunction with eigenvalue $\lambda$. Let $\psi$ be its normal derivative at the boundary. Scaling considerations lead one to expect that the $L^2$ norm of $\psi$ will grow as $\lambda^{1/2}$ as $\lambda \to \infty$. We prove an upper bound of the form $\|\psi \|_2^2 \leq C\lambda$ for any Riemannian manifold, and a lower bound $c \lambda \leq \|\psi \|_2^2$ provided that $M$ has no trapped geodesics (see the main Theorem for a precise statement). Here $c$ and $C$ are positive constants that depend on $M$, but not on $\lambda$. The proof of the upper bound is via a Rellich-type estimate and is rather simple, while the lower bound is proved via a positive commutator estimate. 
\end{abstract}

\maketitle

\section{Introduction} Let $M$ be a bounded domain in $\RR^n$ with
smooth boundary, or more generally a smooth compact Riemannian manifold with
 boundary $\pa M = Y$.  Let $H = -\Delta_M$ be minus the Dirichlet Laplacian on $M$. It is
 well known that $H$ has discrete spectrum $0 < \lambda_1 < \lambda_2 \leq
 \lambda_3 \dots \to \infty$. Let $u_j$ be an $L^2$-normalized eigenfunction corresponding to $\lambda_j$, and let $\psi_j$ be the normal derivative
 of $u_j$ at the boundary. 
In this paper we consider the following question: do there exist constants
$c>0$ and $C < \infty$, depending on $M$ but not on $j$, such that 
\begin{equation}
c \lambda_j \leq \| \psi_j \|^2_{L^2(Y)} \leq C \lambda_j \ \ ?
\label{bounds}\end{equation}
Note that these inequalities exhibit the expected scaling; the square of the
$L^2$ norm of $\psi_j$ involves two derivatives in space, which should
scale as $\lambda_j$ in the high energy limit $\lambda_j \to \infty$. This question was posed by
Ozawa in \cite{Ozawa}. Using heat kernel techniques, he showed
that an averaged version of \eqref{bounds} holds. 
More precisely, he showed that
$$
\sum_{\lambda_j < \lambda} \psi_j^2(y) = \frac{\lambda^{(n/2)+1} }{ (4\pi)^{n/2} \Gamma((n/2)+2)}+ o(\lambda^{(n/2)+1}), \quad \text{for each } y \in Y.
$$
This asymptotic formula (after integrating over $Y$) would be implied by \eqref{bounds} in view of Weyl asymptotics
for the $\lambda_j$. 

In this paper, we prove the upper bound from \eqref{bounds} in general, and
obtain a result for lower bounds. The main theorem is given later in the
introduction. Before stating it, we give some simple examples which show in
particular that the lower bound fails for some Riemannian manfolds with
boundary. Additional, more elaborate examples are given in Section~\ref{Further}.

\

{\it Example 1 --- the disc.} Let $M = \{ x \in \RR^2 \mid |x| < a
\}$ for some $a > 0$. In this case we have an equality
\begin{equation}
\int_{S_1} \psi_j^2(\theta) \, d\theta = \frac{2\lambda_j}{a}.
\label{disc}\end{equation}
This follows from Rellich's identity \eqref{lower-Euc}, which also gives the same result in all dimensions. 

\

{\it Example 2 --- the rectangle.} Let $M = [0,a] \times [0,b]$, where
$a \leq b$. Then normalized eigenfunctions take the form
$$
u = \sqrt{\frac{4}{ab}} \sin \frac{m\pi x}{a} \sin \frac{n\pi y}{b},
\quad m,n \text{ integers}, \quad \lambda = \big( \frac{m\pi}{a} \big)^2 + 
\big( \frac{n\pi}{b} \big)^2.
$$
The square of the $L^2$ norm of its normal derivative is
\begin{equation*}\begin{gathered}
\frac{8}{ab} \Bigg( \int_0^a (\frac{n\pi}{b} \sin \frac{m\pi x}{a})^2
dx  + \int_0^b (\frac{m\pi}{a} \sin \frac{n\pi y}{b})^2 dy \Bigg) \\
= \frac{4}{ab} \Big( \big(\frac{n\pi}{b}\big)^2 a +
\big(\frac{m\pi}{a}\big)^2 b \Big).
\end{gathered}\end{equation*}
From this we see that
\begin{equation}
\frac{4}{b} \lambda \leq \| \psi \|_2^2 \leq \frac{4}{a} \lambda,
\end{equation}
and by taking $m$ or $n$ equal to one and sending the other integer to
infinity, we see that these bounds are the best possible. Thus, for a
square, we get equality as we do for a disc, but in general, the ratio
between the best upper and lower bounds is the ratio of the side lengths. 

\

{\it Example 3 --- the cylinder.} Let $M = [0,a] \times S^1_b$, the product
of an interval with a circle of length $b$. Then the normalized 
eigenfunctions are
$$
u = \sqrt{\frac{2}{ab}} \sin \frac{m\pi x}{a} e^{in\pi\theta/b}, \quad
m,n \text{ integers}, \quad \lambda = \big( \frac{m\pi}{a} \big)^2 + 
\big( \frac{n\pi}{b} \big)^2.
$$
The square of the $L^2$ norm of its normal derivative is
\begin{equation*}\begin{gathered}
2 \ \frac{2}{ab}  \int_0^b \big| \frac{m\pi}{a} e^{in\pi\theta/b} \big|^2 d\theta \\
= \frac{4}{a} (\frac{m\pi}{a})^2.
\end{gathered}\end{equation*}
In this case, the upper bound
$$
\| \psi \|_2^2 \leq \frac{4}{a} \lambda
$$
holds, but no lower bound holds since we may hold $m$ fixed and send $n$ to
infinity. 

\

{\it Example 4 --- the hemisphere.} Let $M$ be the hemisphere
$$
M = \{ x \in \RR^3 \mid |x| = 1, \ x \cdot (0,0,1) \geq 0 \}.
$$
In this case, the eigenfunctions are given by those spherical harmonics
which are odd under reflection in the $(x_1, x_2)$ plane, namely, spherical
harmonics 
$$
u = Y_{lm} = e^{im\phi} P_{lm}(\cos \theta), \quad \lambda = l(l+1),
$$
where $-l \leq m \leq l$ and $l-m$ is odd. We use spherical polar
coordinates $(\theta, \phi)$ where 
\begin{equation}
(x_1,x_2,x_3) =
(r\sin\theta\cos\phi,r\sin\theta\sin\phi,r\cos\theta).
\label{sph-coords}\end{equation}
Let us consider the
case when $m=l-1$. Then 
$$
u = c_l e^{i(l-1)\phi} (\sin\theta)^{l-1} \cos \theta,
$$
and $c$ is a normalization factor. We have 
$$
\| u \|_2^2 = c_l^2 \int_0^\pi \big((\sin \theta)^{l-1} \cos \theta\big)^2
\sin \theta \, d\theta = c_l^2\int_0^\pi \Big( (\sin\theta)^{l-1} -
(\sin\theta)^{l+1} \big) \, d\theta.
$$
A computation shows that 
$$
\int_0^\pi (\sin\theta)^{l+1} \, d\theta = \frac{l}{l+1} \int_0^\pi
(\sin\theta)^{l-1} \, d\theta .
$$
Hence $L^2$ normalization requires that
$$
c_l^{-2} = \frac1{l+1} \frac{l-2}{l-1} \frac{l-4}{l-3} \dots \frac{2}{1} =
\frac{1}{l+1} 2^{-(l+1)} { l-1 \choose (l-1)/2},
$$
where for simplicity we took $l$ to be odd. This has asymptotic
$$
c_l^{-2} \sim l^{-3/2}.
$$
Thus, 
$$
\| \psi \|_2^2 = 4\pi c_l^2 \sim l^{3/2} \sim \lambda^{3/4},
$$
for this class of eigenfunctions. Hence there is no nontrivial lower bound of the form \eqref{bounds} for the
hemisphere. 

\

In all these examples, the upper bound holds, but the lower bound fails in
the last two examples. It is rather clear why the lower bound fails in
example 3; if we hold $m$ fixed and send $n$ to infinity, then in some
intuitive sense the energy in the eigenfunction is concentrating along
lines where $x$ is constant and these never reach the boundary, so most of
the energy is going undetected at the boundary. More precisely, one can
show that the semiclassical wavefront set of the family $u_{m,n}$, with $m$
fixed and $n \to \infty$, is concentrating on geodesics where $x$ is held
fixed. A similar phenomenon is happening for the
hemisphere, although here there is only one geodesic --- the boundary curve
--- which is trapped, so this is a borderline case. Correspondingly, the
lower bound is not violated as strongly in this example. 

These examples lead one to expect that the failure of the lower bound is
related to the presence of geodesics in $M$ which do not reach the
boundary. In fact, we have the following result. 

\begin{thm}\label{main} Let $M$ be a smooth compact Riemannian manifold with boundary. Then the upper bound in \eqref{bounds} holds for some $C$ independent of $j$. The lower bound in \eqref{bounds} holds provided that $M$ can be embedded in the interior of 
a compact manifold with boundary, $N$, of the same dimension,
such that every geodesic in $M$ eventually meets the boundary of $N$. In
particular, the lower bound holds if $M$ is a subdomain of Euclidean
space. 
\end{thm}

We would like to acknowledge helpful conversations with Alan McIntosh, Steve Zelditch and Johannes Sj\"{o}strand.

\section{Rellich-type estimates}\label{sec:Rellich}
To prove the upper bound, and the lower bound for Euclidean domains, we use the following Lemma which we call a Rellich-type estimate.

\begin{lem} Let $u$ be a Dirichlet eigenfunction of $H$. Then for any differential operator $A$, 
\begin{equation}
\int\limits_{M} \langle u, [H,A]u \rangle \, dg = \int\limits_{Y}\frac{\pa u}{\pa \nu} Au \, d\sigma.
\label{Rellich}\end{equation}
\end{lem}

\begin{proof} The proof is very simple. Let $\lambda$ be the eigenvalue corresponding to $u$. We write $[H,A] = [H-\lambda,A]$ and use the fact that $(H-\lambda)u = 0$ to write the integral over $M$ as
\begin{equation*}
\int\limits_{M} \langle (H-\lambda)u, Au \rangle - \langle u, (H-\lambda)Au \rangle \, dg .
\end{equation*}
Then we apply Green's formula; since $u$ vanishes at the boundary, only one of the boundary terms is nonzero, giving \eqref{Rellich}. 
\end{proof}

Let $u$ be an $L^2$-normalized Dirichlet eigenfunction on $M$ with eigenvalue $\lambda$. 
To prove an upper bound for the $L^2$ norm of $\psi = \pa_\nu u$, we choose an operator $A$ so that the right hand side of \eqref{Rellich} is a positive form in $\psi$. To do this, let us choose coordinates $(r,y)$ locally near the boundary such that $r$ is the distance to the boundary; this is a smooth function for $r \in [0, \delta]$, where $\delta$ is sufficiently small. Then we choose $A = \chi(r) \pa_r$, where $\chi \in C_c^\infty(\mathbb{R})$ is identically $1$ for $r$ close to zero, and vanishes for $r \geq \delta$. The right hand side of \eqref{Rellich} is then precisely the square of the $L^2$  norm of $\psi$.
The left hand side may be written (after one  integration by parts)
$$
\int\limits_{M} \langle B_1 u, B_2 u \rangle \, dg 
$$
where $B_i$ are first order (vector-valued) differential operators with smooth coefficients. This may be bounded by
$$
C \int\limits_{M} \langle \nabla u, \nabla u \rangle \, dg = C \lambda,
$$
where $C$ depends on the domain, but not on $\lambda$. This proves the upper bound for any compact Riemannian manifold with boundary. Indeed, we can prove more. Choosing $A$ to be of the form $Q^* Q\pa_r$ near the boundary, where $Q$ is an elliptic differential operator of order $k$ in the $y$ variables, we obtain the upper bound
\begin{equation}
\| \psi \|_{H^k(Y)}^2 \leq C \lambda^{k+1}
\label{deriv-upper-bounds}\end{equation}
for any integer $k$, and hence (by interpolation) any real $k$. 

{\it Remark.} Notice that this argument gives a rather explicit upper bound for $\lim \sup \lambda_j^{-1}  \| \psi_j \|_2$, namely, the reciprocal of the largest $\delta$ for which $r$ is a smooth function on $[0, \delta]$. 

\

To prove the lower bound for Euclidean domains, $M \subset \mathbb{R}^n$, we choose $A$ so that the left hand side, rather than the right hand side, of \eqref{Rellich} is a positive form. We choose 
\begin{equation}
A = \sum_{i=1}^n x_i \frac{\pa}{\pa x_i}.
\end{equation}
As is very well known in scattering theory, the commutator of this with $H$ (which is minus the Euclidean Laplacian here) is $[H,A] = 2H$. Hence, in this case \eqref{Rellich} gives us
\begin{equation}
2\lambda = \int\limits_{Y} \frac{\pa u}{\pa \nu} Au \, d\sigma =
\int\limits_{Y} \nu \cdot x \, \big( \frac{\pa u}{\pa \nu} \big)^2 \, d\sigma
\leq C \| \psi \|_2^2,
\label{lower-Euc}\end{equation}
which gives the lower bound. The equality in \eqref{lower-Euc} was proved by Rellich \cite{Rellich}.

\section{Estimates for eigenfunctions near the boundary}\label{Estimates}
Let $M$ be a smooth compact Riemannian manifold with smooth boundary. 
We choose a coordinate system $(r,y)$ near the boundary where $r$ is
distance to the boundary (a smooth function for small $r$, say $r \leq \delta$) and $y$ is constant on geodesics which hit the boundary normally. In terms of these coordinates, we may write the metric
$$
g = dr^2 + h_{ij} dy_i dy_j
$$
and the Riemannian measure
\begin{equation}
dg = k^2 dr dy, \text{ where } k^4 = \det h_{ij} . 
\label{k}\end{equation}
Let us denote the boundary of $M$ by $Y$, and write $Y_r$ for the set of points at distance $r$ from the boundary, which is a submanifold for $r \leq \delta$. 

Suppose that $u$ is a Dirichlet eigenfunction for $H$, with eigenvalue $\lambda$. 
It will be convenient to change to the function $v = ku$ (this is equivalent to looking at the Laplacian acting on half-densities). Then $v$ solves the equation
\begin{equation}
\pa_r^2 v + \pa_i (h^{ij} \pa_j v) + \lambda v + f v = 0, \quad  h^{ij} = (h_{ij})^{-1}
\label{v-eqn}\end{equation}
where
$$
f = - k^{-1}\pa_r^2 k - k^{-1} \pa_i ( h^{ij} \pa_j k )
$$
is a smooth function on $M$. 

Given an eigenfunction $u$, we define a sort of `energy' $E(r)$ for each value of $r$. This is obtained formally from
the energy for hyperbolic operators, with $r$ playing the role of a time
variable, by switching the sign of the term involving tangential
derivatives. Let
\begin{equation}
E(r) = \frac1{2} \int_{Y_{r}} \big(  v_r^2 + (\lambda + f)v^2 -
h^{ij} \pa_i v \pa_j v \big) dy.
\end{equation}

\begin{lem} For $r \in [0, \delta]$, we have an estimate
\begin{equation}
|E(r)| \leq C\lambda
\label{energy-est}\end{equation}
where $C$ is independent of $\lambda$.
\end{lem}

\begin{proof}
From the upper bound argument in the previous section, we know that $E(0) = \frac1{2} \| \psi \|^2 \leq C\lambda$. We compute the derivative of $E(r)$:
\begin{equation}\begin{gathered}
\frac{\pa}{\pa r} E(r) = \int_{Y_{r}} \big( \pa_r^2 v \pa_r v + 
(\lambda + f)v \pa_r v  -
h^{ij} \pa_i v \pa_r \pa_j v \\ 
+ \frac{\pa h^{ij}}{\pa r} \pa_i v
\pa_j v + \frac{\pa f}{\pa r} v^2 \big) dy.
\end{gathered}\end{equation}
Integrating by parts in the third term, and using the equation for $v$, we
obtain 
\begin{equation}
\big| \frac{\pa E}{\pa r} \big|(r) \leq C \int_{Y_{r}} \big( v^2 + |\nabla
v|^2 \big) dy \leq C \int_{Y_{r}} \big( u^2 + |\nabla u|^2 \big)
k^2 dy.
\end{equation}
Thus, for $r_0 \in [0, \delta]$, 
$$
E(r_0) = E(0) + \int_0^{r_0} \frac{d}{dr} E(r) dr \ \leq \ C\lambda + \int_M \big( u^2 + |\nabla u|^2 \big) dg \ \leq \ C\lambda.
$$
\end{proof}

Next, we derive an estimate on the $L^2$ norm of $u$ on $Y_r$, exploiting the fact that $u$ vanishes on the boundary. 

\begin{lem}  There exists $C > 0$, independent of $\lambda$, such that
\begin{equation}
\int_{Y_r} u^2 d\sigma(y) \leq C\lambda r^2 \text{ for all } r \in [0, \frac{\delta}{3}].
\label{bdy-est}\end{equation}
\end{lem}

{\it Remark.} Given $u$, this bound follows, for small $r$, from Taylor's theorem and the upper bound from Section~\ref{sec:Rellich}. The point of this Lemma is that the estimate holds for $r$ in a fixed interval, {\it independent of $\lambda$}.

\begin{proof}
Consider the $L^2$ norm on $Y_r$:
\begin{equation}
L(r) = \int_{Y_r} u^2 k^2 dy = \int_{Y_r} v^2 dy.
\end{equation}
Then we have
$$ L'(r_0) = \int_{Y_{r_0}} 2 v v_r\ dy$$
and
\begin{align*}
 L''(r_0) &= 2 \int_{Y_{r_0}} \big( v_r^2 + v v_{rr} \big) \ dy \\
&= 2 \int_{Y_{r_0}} \big( v_r^2 - v \partial_i (h^{ij} v_j) - v \lambda v - v f v\big) \ dy \\
&= 2 \int_{Y_{r_0}} \big( v_r^2 + h^{ij} v_i v_j -  \lambda v^2 - f v^2 \big) \ dy \\
&= 4 \int_{Y_{r_0}}  v_r^2  \ dy \ -  4E(r_0).
\end{align*}
On the other hand, from Cauchy-Schwarz we have
$$ 4 \int_{Y_{r_0}} v_r^2 \ \geq \ \frac{(\int_{Y_{r_0}} 2 v v_r\ dy)^2}{\int_{r=r_0} v^2\ dy} = \frac{L'(r_0)^2}{L(r_0)}.$$
Thus we have the differential inequality for $L(r)$:
\begin{equation}
L'' \geq \frac{(L')^2}{L} - C \lambda ,
\label{diff-ineq}\end{equation}
for some constant $C$ depending only on the manifold $M$.

Intuitively, this inequality says that if $L'$ ever becomes too big, then $L''$ will also become big.  This should create a feedback loop which causes $L$ to grow rapidly, which we know cannot happen because we have 
\begin{equation}
\int_0^\delta L(r)  dr \leq  \int_M u^2 \, dg = 1. 
\label{u-bd}\end{equation}
So we should be able to get upper bounds on $L'$ from \eqref{diff-ineq} and \eqref{u-bd}, which in turn should imply an upper bound on $L$.

We use an integrating factor.  After some rearrangement, \eqref{diff-ineq} becomes
$$
\frac{2 L' L''}{L^2} - \frac{2 (L')^3}{L^3} + \frac{2 C\lambda L'}{L^2} \geq 0.$$
Since the left-hand side is the derivative of the quantity
$$
B(r) := \frac{L'(r)^2}{L(r)^2} - \frac{2C \lambda}{L(r)}$$
we thus see that $B$ is non-decreasing.

Suppose that $B(r_0)$ and $L'(r_0)$ were positive for some $0 < r_0 < \delta/3$.  Then we would have $B(r) > 0$ for all $r \geq r_0$, so
$$ L'(r)^2 > 2C\lambda L(r) \hbox{ for all } r \geq r_0.$$
In particular $L'(r)$ would be strictly positive for $r \geq r_0$.
We rearrange this as
$$ (L(r)^{1/2})' = \frac{1}{2} L'(r) L(r)^{-1/2} > \sqrt{\frac{C}{2}\lambda} \hbox{ for all } r > r_0.$$
This would give
$$ L(r) \geq C' \lambda \hbox{ for all } r \geq \frac{2\delta}{3}.$$
This would contradict the bound \eqref{u-bd} for large $\lambda$. Thus, for $\lambda \geq \lambda_0$ (where $\lambda_0$ depends only on $M$), we must have $B(r) \leq 0$ or $L'(r) \leq 0$ for all $0 < r  \leq \delta/3$. In either case
$$ (L(r)^{1/2})' = \frac{1}{2} L'(r) L(r)^{-1/2} \leq \sqrt{\frac{C}{2}\lambda} \hbox{ for all } r \leq \frac{\delta}{3}, \quad \lambda \geq \lambda_0.
$$
Since $L(0) = 0$, this implies \eqref{bdy-est}. 
\end{proof}

\section{The lower bound on Riemannian manifolds}
To find a lower bound on an arbitrary compact Riemannian manifold, $M$ satisfying the conditions of the theorem, we need to find an operator which has a positive commutator with $H$. We begin by constructing a first order pseudodifferential operator $A$ on $N$ which has this property to leading order, ie such that the symbol of $i[H,A]$ is positive. 

\begin{lem}\label{Q-lemma} Given any geodesic $\gamma$ in $S^*N$, there is a first order, classical, self-adjoint pseudodifferential operator $Q$ satisfying the transmission condition (see \cite{Ho}, section 18.2), and properly supported on $N$, such that the principal symbol $\sigma(i[H,Q])$ of $i[H,Q]$ is nonnegative on $T^*M$, and 
\begin{equation}
\sigma(i[H,Q]) \geq \sigma(H) = |\xi|^2 \ \text{ on a conic neighbourhood $U_\gamma$ of } \gamma \cap T^*M.
\label{comm-cond}\end{equation}
\end{lem}

\begin{proof}
Recall that the principal symbol of $i[H,Q]$ is given by the Hamilton vector field of $H$ applied to $q = \sigma(Q)$, the principal symbol of $Q$. 
In the case of $H =
-\Lap_g$, this is $|\xi|$ times the arc-length derivative along
geodesics, where $|\xi| = |\xi|_g = \sqrt{g^{ij}\xi_i \xi_j}$; note that
$|\xi|$ is constant under geodesic flow. 

To construct $q$, choose a point $p$ on $\gamma$ and a small piece of
submanifold $V$ transversal, within $S^*N$,  to $\gamma$. 
We next choose functions $\psi$ and $b$ which are homogeneous of
degree zero, so we only need to specify them on $S^*N$. 
We first specify $\psi$ on $V$, requiring it to be identically $1$ near $p$
and zero outside a slightly larger neighbourhood of $p$. We then 
extend $\psi$ to $S^*N$ be requiring it to be constant along geodesics
intersecting $V$, and zero on any geodesic not intersecting $V$. 
Thus $\psi$ is supported in a small conic neighbourhood $O$ of $\gamma$. 
We define $b$ on the support of $\psi$ in $S^*N$ by prescribing it to be
zero on $V$ and to have derivative equal to $1$ along geodesics, with
respect to arc length. Then we define $a$, a function homogeneous of degree
zero on $T^*N$, by $a = \psi b$. Due to our geodesic condition on $M$, $a$
is a smooth function on $T^*N \setminus 0$. 

We next choose cutoff functions as follows: let
$\chi$ be a smooth
function on $N$ which is $1$ on $M$ and is supported in the interior of
$N$, and let $\chi_2$ be a smooth function on $T^*N$ which is $1$ when
$|\xi| \geq 1$ and $0$ when $|\xi| \leq 1/2$. 
We define symbols $q'$ and $q''$ on $T^*N$ by $q' = |\xi| \chi \chi_2
a$, and $q'' = -q'(x,-\xi)$. Then $q''$ is supported in a small conic neighbourhood of $-\gamma$ and its symbol is increasing along $-\gamma$. Let $q = q' + q''$; since $q$ is odd, we can find a classical pseudodifferential operator $Q$ with symbol $q$ satisfying the transmission condition. Replacing $Q$ by $\frac1{2} (Q + Q^*)$ if necessary, we may assume that $Q$ is self-adjoint. By construction, $Q$ satisfies \eqref{comm-cond}. 
\end{proof}

We now use Lemma~\ref{Q-lemma} to construct our operator $A$. For each geodesic $\gamma$ in $S^*N$, we have a conic neighbourhood $U_\gamma$ as in the Lemma. By compactness of $S^*M$, a finite number of the $U_\gamma$ cover $S^*M$. Let $A$ be the sum of the corresponding $Q_\gamma$. Then Lemma~\ref{Q-lemma} implies that 
\begin{equation}
\sigma(i[H,A]) \geq |\xi|^2 \ \text{ on } T^*M . 
\label{A}\end{equation}

Our next task is to get an identity for $A$ similar to \eqref{Rellich}. Since $A$ is now \discretionary{pseudo-}{differential}{pseudodifferential} and therefore non-local, it is important to take into account the fact that $u$ is defined only on $M$. Let $\ubar$ denote the extension by zero to $N$. 

\begin{lem} Let $u$ be a Dirichlet eigenfunction for $H$, and let $A$ be a first order classical pseudodifferential operator satisfying the transmission condition. Then
\begin{equation}
\int\limits_{M} \langle \ubar, [H,A] \ubar \rangle \, dg = 2 \Imag \int\limits_{Y} \frac{\pa u}{\pa \nu} A \ubar \, d\sigma - \int\limits_{Y} \big( \frac{\pa u}{\pa \nu} \big)^2 c \, d\sigma,
\label{pseudo-Rellich}\end{equation}
where $c(y) = \lim_{\rho \to \infty} \rho^{-1} a(0,y,\rho,0)$. \end{lem}

{\it Remark. } In \eqref{pseudo-Rellich}, the restriction of $A \ubar$ to $\pa M$ is taken from the interior of $M$; this is well defined since $A$ satisfies the transmission condition and $u$ is smooth on $\overline{M}$ (see \cite{Ho}).

\begin{proof}
Note that $\ubar \in H^{3/2-\epsilon}(N)$ for every $\epsilon > 0$, but it is not in $H^{3/2}$. Consequently, we cannot unwrap the commutator $[H,A]$ directly --- since $HA$ and $AH$ are third order operators --- but must use a limiting argument. Let $\phi(t)$ be a smooth function which is zero for $t < 1$ and one for $t > 2$. Choose coordinates $(r,y)$ as above. Notice that for $0 < \epsilon < \delta/2$, $\phiep(r) = \phi(r/\ep)$ may be regarded as a smooth function on $N$, such that $\uep = \ubar \phiep$ is a smooth  function on $N$. 

Since $[H,A]$ is second order, and $\uep \to \ubar$ in $H^1(N)$ as $\ep \to 0$, we have 
$$
\int_{M} \langle \ubar, [H,A] \ubar \rangle \, dg = \lim_{\ep \to 0} 
\int_{M} \langle \uep, [H - \lambda,A]\uep \rangle \, dg.
$$
Now we can unwrap the commutator, getting 
$$
\int_{M} \langle (H-\lambda)\uep, A \uep \rangle - \langle A \uep, (H - \lambda) \uep \rangle \, dg = 2 \Imag \int_{M} \langle (H-\lambda)\uep, A \uep \rangle \, dg.
$$
We write 
\begin{equation}
\langle (H-\lambda)\uep, A \uep \rangle = \langle (H-\lambda)\uep, A \ubar \rangle -
\langle (H-\lambda)\uep, A ((1-\phiep)u) \rangle.
\label{chopping}\end{equation}
Now we use the fact that $A$ satisfies the transmission condition. This implies that $A\ubar$ is smooth up to the boundary of $M$, ie $A\ubar \restriction \overline{M}$ is smooth. If we compute $(H-\lambda)\uep$, we get
$$
(H-\lambda) \uep = -\pa_r u \pa_r \phiep - u \pa_r^2 \phiep + 2 u k_r \pa_r \phiep.
$$
The distributional limit of this as $\ep \to 0$ is $- \psi \delta_{\pa M}$, and it is always supported in the interior of $M$. Therefore, 
$$
\int_{M} \langle (H-\lambda)\uep, A \ubar \rangle \, dg \to \int_{Y} \psi A \ubar \, d\sigma
$$
where $A \ubar \restriction M$ is taken as the limit from the interior of $M$. 

To deal with the second term in \eqref{chopping}, we decompose $A$. Let $a$ be the principal symbol of $A$, and let $y$ be a point on the boundary of $M$.  Let $(\rho, \eta)$ be cotangent coordinates dual to $(r,y)$. Then, restricted to the line $\eta = 0$, the symbol is the same as that of $c(y) D_r$ for $c$ as defined in the Lemma. Since the symbols $\eta_i$ of $D_{y_i}$ vanish on this line, and have independent differentials there, we can write the symbol $a$ as 
$$
a = c(y) \rho + \eta_i b_i,
$$
for some symbols $b_i$ of order zero. Correspondingly, we get a decomposition
$$
A = c D_r + B_i D_{y_i} + B',
$$
where $B_i$ and $B'$ are pseudodifferential operators of order zero. Then replacing $A$ by $c D_r$ in  the second term in \eqref{chopping} gives 
$$
-\int_{M} \langle (H-\lambda)\uep, c D_r (u (1 - \phiep)) \rangle \, dg \to -i\int_{Y} c(y) \psi^2 d\sigma.
$$
Thus, to prove the Lemma it remains to show that 
$$
\int_{M} \langle (H-\lambda)\uep, (B_i D_{y_i} + B') (\ubar (1 - \phiep)) \rangle \, dg \to 0.
$$
Notice that, as $\ep \to 0$, the $L^2$ norm of $(H-\lambda) \uep$ is $O(\ep^{-1/2})$, while the $L^2$ norms of $D_{y_i} (\ubar ( 1 - \phiep))$ and of $ \ubar( 1 - \phiep)$ are $O(\ep^{3/2})$. Since zeroth order pseudodifferential operators are bounded on $L^2$, this term indeed goes to zero. This completes the proof of the Lemma.
\end{proof}

Now we use \eqref{pseudo-Rellich}  and \eqref{bdy-est} to prove the lower bound.
First, we deal with the left hand side of \eqref{pseudo-Rellich}. 
Using \eqref{A}, we may write 
$$
i[H,A] = H + B^* B + R_1 + R_2,
$$
where $B$ and $R_1$ are first order pseudodifferential operators, and $R_2$ is second order with kernel smooth on $M \times M$. Thus, we can estimate from below
\begin{equation}\begin{gathered}
\int_{M} \langle \ubar, i[H,A] \ubar \rangle \, dg  \geq \lambda + 
\int_{M} \langle \ubar, (R_1 + R_2) \ubar \rangle \, dg \\
\geq \lambda - C \| u \|_{H_1} \| u \|_{L^2} = \lambda - C \lambda^{1/2}.
\end{gathered}\end{equation}
We want to show that this is no bigger than $C \| \psi \|_2^2$. To do this, it is sufficient to show that 
\begin{equation}
\| A \ubar _{\restriction Y}  \|_{L^2(Y)}^2 \leq C \lambda + C \| \psi \|^2,
\label{suff}\end{equation}
since by Cauchy-Schwarz we can then estimate the right hand side of \eqref{pseudo-Rellich} by
$$
\int\limits_{Y} \langle \psi, A \ubar \rangle \, d\sigma \leq \frac{1}{\ep} \| \psi \|_{L^2}^2 + \ep \| A \ubar_{ \restriction Y }\|_{L^2(Y)}^2,
$$
use \eqref{suff} and absorb the $\| Au \|^2$ term in the left hand side. This, to finish the proof of the lower bound, it is enough to prove \eqref{suff}. 

\

To do this, we require the following Lemma, which is due to Boutet de Monvel \cite{BdM} and Vi\v{s}ik  and E\v{s}kin \cite{VE}. 

\begin{lem}\label{BdM} Let $B$ be a classical pseudodifferential operator on $N$ of integral order $k \geq -1$ satisfying the transmission condition. Let $B_r$ denote the operator from $L^2(Y_r)$ to $L^2(Y)$ given by $B_r(f_r) = B(f_r \otimes \delta_{Y_r}) \restriction Y$. (For $r=0$, we take 
$\lim_{\ep \downarrow 0} (B(f_0 \otimes \delta_{Y}))(\ep, \cdot)$; this is well defined for operators satisfying the transmission condition.) Then we have an estimate
\begin{equation}
\| B_r \|_{L^2(Y_r) \to L^2(Y)} \leq C r^{-k-1}, \ r \in (0, \delta].
\label{trans-est}\end{equation}
For $k=-1$ we have $B_r$ uniformly $L^2$ bounded down to and including $r=0$; indeed, in this case, if we identify $Y_r$ and $Y$,  then $\{ B_r \}$ is a uniformly bounded family of zeroth order pseudodifferential operators on $Y$, for $r \in [0, \delta]$. 
\end{lem}

\begin{proof} The proof of this is at least implicitly contained in \cite{Ho}, in the discussion before Theorem 18.2.17, but we give a self-contained proof here which is based on that in \cite{Ho}. First, since $\Psi^k(N)$ is invariant under adjoints, it is sufficient to bound $B$ mapping from $L^2(Y)$ to $L^2(Y_r)$, which by an abuse of notation we shall denote $B_r$ for the duration of this proof. 

Let the kernel of $B$ be given as an oscillatory integral near $r = r' = 0$ by
$$
B(r,y,r',y') = \int e^{i(r-r')\xi_1} e^{i(y-y') \cdot \xi'} b(r,y,\xi_1,\xi') \, d\xi_1 \, d\xi'.
$$
The kernel of $B_\ep$ is obtained by restricting this to $r=\ep$, $r' = 0$. 
Then we wish to obtain a uniform $L^2$ bound on $r^{k+1} B_r$. Multiplying the kernel by $r^{k+1}$ is the same as applying $D_{\xi_1}^{k+1}$ to the phase. Integrating these derivatives by parts reduces the order of the symbol to $-1$. Thus, it is enough to consider the case $k=-1$. 

Now consider the operator $B_r$ acting on a fixed function $v \in L^2(Y)$. Equivalently, we can look at $B$ acting on $v \delta_{Y}$. This is the distributional limit of $B$ acting on $v \ep^{-1} \phi(r/\ep)$, where $\phi \in C_c^\infty(\RR)$ has integral one and is supported in $|r| \leq 1$, say.  The result is
\begin{equation}
\int e^{ir\xi_1} e^{iy\cdot\xi'} b(r,y,\xi) \hat v(\xi') \hat \phi(\ep \xi_1) \, d\xi_1 \, d\xi'.
\label{eq2}\end{equation}
Thus, the symbol of $B_r$ is 
\begin{equation}
\lim_{\ep \to 0} \int e^{ir\xi_1} b(r,y,\xi_1, \xi') \hat \phi(\ep \xi_1) \, d\xi_1 .
\label{symbol}\end{equation}
If $b$ in \eqref{symbol} were a symbol of order $-2$, then the limit in \eqref{symbol} is
\begin{equation}
\int e^{ir\xi_1} b(r,y,\xi_1, \xi')  \, d\xi_1 .
\end{equation}
which is a symbol of order $-1$ in the $\xi'$ variables, depending uniformly on $r$ (ie, with uniform estimates on all seminorms). For $b$ of order $-1$, we choose a cutoff function $\chi(\xi_1)$,  where $\chi$ is supported in $|\xi_1| \leq 1$ and is identically one near $\xi_1 = 0$, and write $b = \chi(\xi_1) b + b'$. Then $\chi(\xi_1) b$ is a symbol of order $-1$ in the $\xi'$ variables, uniformly in $r$ and $\xi_1$ (since it is compactly supported in $\xi_1$), so this integral gives a symbol of order $-1$ uniformly in $r$. Finally, $b'$ is supported away from $\xi_1 = 0$, so we may write using Taylor's theorem
$$
b' = (1 - \chi(\xi_1))\Big( c(r,y) \xi_1^{-1}  + \sum_{j=2}^n \xi_j' b_j \Big),
$$
where $b_j$ are symbols of order $-2$. Using the statement above about symbols of order $-2$, we see that the $\xi_j' b_j$ terms give symbols of order $0$, uniformly in $r$. On the other hand, the $\xi_1^{-1}$ term can be calculated explicitly. In the integral \eqref{symbol}, with $b$ replaced by $(1 - \chi(\xi_1)) c(r,y) \xi_1^{-1}$, the integral can be moved into the upper half complex $\xi_1$ plane provided $\ep < r$, and then the result is 
$$
c(r,y) \Big( \int_{-1}^1 e^{ir\xi_1} (1 - \chi(\xi_1)) \xi_1^{-1} d\xi_1 \ \  + \!\!\!\!\!\! \int\limits_{|\xi_1| = 1 ,\  0 < \arg \xi_1 < \pi}  \!\!\!\!\!\! e^{ir\xi_1} \xi_1^{-1} d\xi_1 \Big),
$$
which is a multiplication operator depending smoothly on $r$. This completes the proof. 
\end{proof}

Let $H^{-1}$ denote the inverse of the Laplacian, with Dirichlet boundary conditions, on $N$. Then, in a neighbourhood of $M$, $H^{-1}$ is a classical pseudodifferential operator of order $-2$ satisfying the transmission condition. Since $H \ubar = \lambda \ubar + \psi \delta_{Y}$, we find $\ubar = H^{-1} ( \lambda \ubar + \psi \delta_{Y})$. Thus, 
$$
A\ubar = AH^{-1} (\lambda \ubar + \psi \delta_{Y}).
$$
By Lemma~\ref{BdM},  $AH^{-1}(\psi \delta_{Y}) \restriction Y$ is bounded by
$C \| \psi \|_2$, which satisfies \eqref{suff}. Let us consider $\lambda AH^{-1} \ubar $.  
To do this we break up $\ubar$ into a `close' and `far' part. Since the length scale is $\lambda^{-1/2}$, let us use the functions $\phiep$ and $\uep$ as above, but choosing the value $\ep = \lambda^{-1/2}$. Consider first $\lambda AH^{-1} (\ubar (1 - \phiep)) $.
By Lemma~\ref{BdM}, we get an estimate independent of $r$ on the operator norm of $(AH^{-1})_r$, and $\ubar(1-\phiep)$ is supported in $r < 2\lambda^{-1/2}$, so we can estimate using \eqref{bdy-est}
$$
\| \lambda AH^{-1} (u(1-\phiep)) \|_2 \leq C \lambda \int_0^{2\lambda^{-1/2} } r \lambda^{1/2} dr \leq C \lambda^{1/2}, 
$$
which satisfies \eqref{suff}. Finally, consider the term $\lambda AH^{-1} \uep$. Writing $\uep = (\lambda^{-2} H^2 u)\phiep$, we obtain
$$
\lambda AH^{-1} \uep = 
\lambda^{-1} AH \uep - \lambda^{-1} A H^{-1} [H^2, \phiep] u.
$$
The first term is estimated via Lemma~\ref{BdM} and \eqref{bdy-est}, and using smoothness of the kernel of $AH$ at $(r,y,r',y')$ for $r=0$ and $r' \geq \delta/3$, by
$$
\lambda^{-1} \Big(  \int_{\lambda^{-1/2}}^{\delta /3} r^{-4} \lambda^{1/2} r dr + O(1) \Big) \leq C \big( \lambda^{1/2} + \lambda^{-1} \big)
$$
which satisfies \eqref{suff} for large $\lambda$. The terms coming from the commutator are treated similarly. 

This completes the proof of Theorem~\ref{main}.

\section{An alternative proof of the lower bound}
There is an alternative way to prove the lower bound, which depends on a trick
due to Morawetz, Ralston and Strauss \cite{MRS}. First, we recall some facts about spherical harmonics. Let $\Delta_{S^{n-1}}$ denote the Laplacian on the $(n-1)$-sphere. Then $-\Delta_{S^{n-1}}$ has eigenvalues $k(n+k-2)$, $k = 0, 1, 2, \dots$. Let the corresponding eigenspace be denoted $V_k$; elements of $V_k$ are called spherical harmonics. We recall that  for every $\phi \in V_k$, the function $r^k \phi$ (thought of as a function on $\RR^n$ written in polar coordinates) is a homogeneous polynomial, of degree $k$, on $\RR^n$. 

The alternative proof goes as follows. 
We use the symbol $a$ of the operator $A$ constructed in the previous section, restrict it to the cosphere bundle of $N$, and expand in spherical harmonics. Note that since $a$ is odd, and spherical harmonics are even or odd with $k$, all the terms with $k$ even vanish. Thus we have
$$
a \restriction S^*N= \sum_{l=0}^\infty  \phi_{2l+1}(x, \frac{\xi}{|\xi|}), \quad \phi_k(x, \cdot) \in V_k(S^*_x N).
$$
Since $a$ is smooth, this sum converges in $C^\infty$, and certainly in $C^1$. Therefore, there is a $K$ such that the operator $A'$ with symbol 
$$
a' = \sum_{l=0}^{K-1} \phi_{2l+1}(x, \frac{\xi}{|\xi|})
$$
also has positive commutator with $H$. However, following \cite{MRS}, we can turn $A'$ into a differential operator $P$ of order $2K-1$, by letting 
$$
p = \sigma(P) = \sum_{l=0}^{K-1} \phi_{2l+1}(x, \frac{\xi}{|\xi|}) |\xi|^{2K- 1}.
$$
The symbol $p$ is a polynomial on each fibre of $T^*N$, so we may take $P$ to be a differential operator. Moreover, the symbol of $i[H,P]$ satisfies
$$
\sigma(i[H,P]) = |\xi|^{2K} \big( \sigma(i[H,P]) |_{|\xi| = 1} \big) \geq c|\xi|^{2K}
\quad \text{ for some } c > 0.
$$

We now apply apply the G{\aa}rding inequality to $Q = i[H,P]$. We need to apply the form which is valid for $u \in H^K(M)$ as opposed to $H^K_0(M)$ which is the standard situation. This gives
\begin{equation}
\int_M \langle u, Qu \rangle dg \geq c \| u \|_{H^K(M)}^2 - C \Big(  \| u \|_{L^2(M)}^2 + \sum_{k=0}^{K-1} \| \pa_r^k u \| _{H^{K-1/2-k}(Y)}^2 \Big),
\label{Garding}\end{equation}
where $c$ is a positive constant depending on $P$ and $(M,g)$. 
This can be derived from the usual G{\aa}rding inequality by writing $u = v + w$, where $v$ solves the Dirichlet problem $(P + \alpha)v = 0$ with the same Dirichlet data as $u$. 
(Here $\alpha$ is a sufficiently large constant.) We then get an estimate on the $H^K$ norm of $v$ in terms of the boundary norms of $u$ appearing in \eqref{Garding}. On the other hand, we can apply the standard G{\aa}rding inequality to $w$. 

Thus, \eqref{Rellich} and \eqref{Garding} combined yield
\begin{equation*}
\| u \|_{H^K(M)}^2 \leq C + C \sum_{k=0}^{K-1} \| \pa_r^k u \| _{H^{K-1/2-k}(Y)}^2
+ \big| \int_{Y} \langle \psi, Pu \rangle \, d\sigma \big|.
\end{equation*}
We can estimate the left hand side from below by a constant times $\lambda^K$. 
Consider the right hand side. Let us write $u = k^{-1}v$, where $k$ is as in \eqref{k}, so that $v$ satisfies \eqref{v-eqn}. Then $Pu = \tilde P v$, where $\tilde P = P \circ k$ is a differential operator of order $2K-1$. Since $k$ is smooth, we obtain 
\begin{equation}
\lambda^K \leq C + C \sum_{k=0}^{K-1} \| \pa_r^k v \| _{H^{K-1/2-k}(Y)}^2
+ C \big| \int_{Y} \langle \psi, \tilde P v \rangle \, d\sigma \big|.
\label{eqqq}\end{equation}
Since $v=0$ at $Y$, and we are interested in $\tilde P v \restriction Y$, we may assume that $\tilde P = P' \pa_r$, where $P'$ has order $2K-2$. Using \eqref{v-eqn}, we may replace $\pa_r^2 v$ by $-(\lambda - f)v - \pa_{y_i} (h^{ij} \pa_{y_j} v)$ repeatedly, until only $\pa_r \pa_y^\alpha v$ terms remain. Thus we have 
\begin{equation*}
\tilde P v \restriction Y = \sum_{j=0}^{K-1} \lambda^j P_j \psi,
\end{equation*}
where $P_j$ is a differential operator on $Y$ of order $2(K-1-j)$, independent of $\lambda$. Hence \eqref{eqqq} becomes
\begin{equation}
\lambda^K \leq C + C \sum_{k=0}^{K-1} \lambda^{k-1} \| \psi \| _{H^{K-1/2-k}(Y)}^2
+ C \sum_{j=0}^{K-1} \lambda^{j}\big| \int_{Y} \langle \psi, P_j \psi \rangle \, d\sigma \big|.
\end{equation}
Using the upper bound estimate \eqref{deriv-upper-bounds} on the sum over $k$ and for all terms in the sum over $j$ with $j < K-1$, we find
\begin{equation}
\lambda^K \leq C + C (1 + \lambda^{K-1/2}) + C \lambda^{K-1} \| \psi \|_{L^2}^2 
+ C \lambda^{K-1/2} \| \psi \|_{L^2} .
\label{almost}\end{equation}
Finally, we estimate the last term in \eqref{almost} by $\ep \lambda^{K}
 + \ep^{-1} C^2 \lambda^{K-1} \| \psi \|_2^2$, and absorb the $\ep \lambda^K$ term in the left hand side. This proves the lower bound.

\section{Further examples}\label{Further}
Here we briefly discuss further examples which illustrate how the presence of trapped geodesics in the interior of $M$ affects the behaviour of eigenfunctions at the boundary. 

\

{\it Example 5 --- the spherical cylinder.} Consider the manifold $M \subset S^2$ given in terms of coordinates $(\theta, \phi)$ as in \eqref{sph-coords}, by
$$
M = \{ (\theta, \phi) \mid \frac{\pi}{2} - a \leq \theta \leq \frac{\pi}{2}  + a, \quad \phi \in S^1 \}, \quad 0 < a < \frac{\pi}{2}.
$$
Let $y = \theta - \pi/2$. Here, the equator at $y=0$ is a trapped geodesic, as is every nearby geodesic, so this is an extreme case where there is a nonempty open set of periodic geodesics. In this case, there is a sequence of normalized Dirichlet eigenfunctions $u_l$ where the lower bound fails spectacularly, namely there is exponential decay of $\| \psi_l \|_2$:
\begin{equation}
\| \psi_l \|_2^2 \leq C e^{-cl}
\label{exp-decay}\end{equation}
This follows by separating variables, $u = e^{il\phi} v(y)$, and reducing to a one dimensional problem where $l^{-1}$ plays the role of a semiclassical parameter. The positive curvature of the sphere manifests itself as a potential function with a nondegenerate minimum at $y=0$. Well known estimates (see e.g. \cite{Helffer}, chapter 3) show that there are normalized eigenfunctions $u_l$ which localize near $y=0$ as $l \to \infty$. More precisely, for any $\delta > 0$ and differential operator $P$ there are constants $c > 0$ and $C < \infty$ such that $| P u_l | \leq C e^{-c l}$ in the set $\{ |y| > \delta \}$. If we then apply our upper bound argument from Section~\ref{sec:Rellich}, which involves only the gradient of the eigenfunction near the boundary, we obtain \eqref{exp-decay}.

\

{\it Example 6 --- a finite hyperbolic cylinder.} Let $M$ be the Riemannian manifold
$$
M = [-a,a]_y \times S^1_\phi , \quad g = dy^2 + (\cosh y)^2  d\phi^2.
$$
Again, there is a periodic geodesic at $y=0$. However,
this is the opposite extreme to the example above. Now $M$ has negative curvature, the geodesic is unstable, and it is in fact the unique trapped geodesic on $M$. (Other geodesics may be trapped as $t \to \infty$ or $t \to -\infty$, but not both.) Using results of Colin de Verdi\`ere and Parisse \cite{CdV-P} one can show that
there is a sequence of normalized eigenfunctions $u_l$, $l = 1, 2, 3, \dots$, such that
\begin{equation}
 \frac{c\lambda}{\log \lambda} \leq \| \psi_l \|_2^2 \leq \frac{C\lambda}{\log  \lambda}.
\label{hyp}\end{equation}
Thus, the lower bound is violated, but just barely. Both this example and the previous one are based on sequences of eigenfunctions which microlocally concentrate at a periodic geodesic. Since the geodesic is in the interior of the manifold, the eigenfunctions are concentrated in the interior, and so are not as large near the boundary as one would normally expect. 

\

{\it Example 7 --- Neumann eigenfunctions.} One might hope to prove bounds of the form
$$
c \leq \| u_{| Y} \|_{L^2(Y)} \leq C
$$
for normalized Neumann eigenfunctions $u$. The following examples show that neither the upper nor lower bound is valid in general. 

Consider the unit disc, as in Example 1. There is a basis of Neumann eigenfunctions of the form 
$$
u_{k,n} = c_{k,n} e^{in\theta} J_n(j'_{k,n}r),
$$
where $c_{k,n}$ is a normalization constant and $j'_{k,n}$ denotes the $k$th zero of the derivative of the $n$th Bessel function $J_n$. The eigenvalue for this eigenfunction is $\lambda = (j'_{k,n})^2$. Denoting the boundary value of this eigenfunction by $\chi_{k,n}$, Rellich's identity \eqref{Rellich} gives us
\begin{equation}
2\lambda = -\int_{S^1} \langle u,  \pa_r^2 u \rangle \, d\theta = \int_{S^1} \langle u, ( \lambda + \pa_\theta^2) u \rangle\, d\theta =(\lambda - n^2) \| \chi_{k,n} \|_2^2.
\label{Neumann-disc}\end{equation}
This shows that there is a lower bound for $\| \chi_{k,n} \|_2$. However, the first positive critical point of $J_n$, $j'_{1,n}$, is known to have asymptotic 
$$
j'_{1,n} \sim n + c n^{1/3}, \quad n \to \infty.
$$
So fixing $k = 1$ and sending $n \to \infty$ in \eqref{Neumann-disc} shows that there is no upper bound for $\| \chi_{k,n} \|_2$.

We remark that it is not hard to derive a lower bound for {\it convex} Euclidean domains from Rellich's identity. 
 On the other hand, the argument of Example 5 applies equally to Neumann eigenfunctions, so the lower bound is not generally valid. It is not clear to the authors whether the lower bound holds under the nontrapping assumption of Theorem~\ref{main}.

\

{\it Example 8 --- vector fields are not enough.} One might wonder whether, on an arbitrary compact Riemannian manifold, with no trapped geodesics (that is, satisfying the conditions of Theorem~\ref{main}), one could choose a first order {\it differential} operator $A$
whose commutator with $H$ had a positive symbol.  In this section, we show via an example that this is impossible in general. First we analyze what it means for a vector field to have a positive commutator with $H$. Let the symbol of $A$ be $a_i(x) \xi_i$. The Hamilton vector field of $H$ is $\xi_i \pa_{x_i}$, so having a positive commutator requires that
\begin{equation}
\xi_i \xi_j \frac{\pa a_j}{\pa x_i} > 0 \text{ for } |\xi| \neq 0.
\label{pos-def}\end{equation}
Thus, the matrix $\pa_{x_i} a_j$ must be positive definite. 

\vskip 20pt
\centerline{\includegraphics[width=2in]{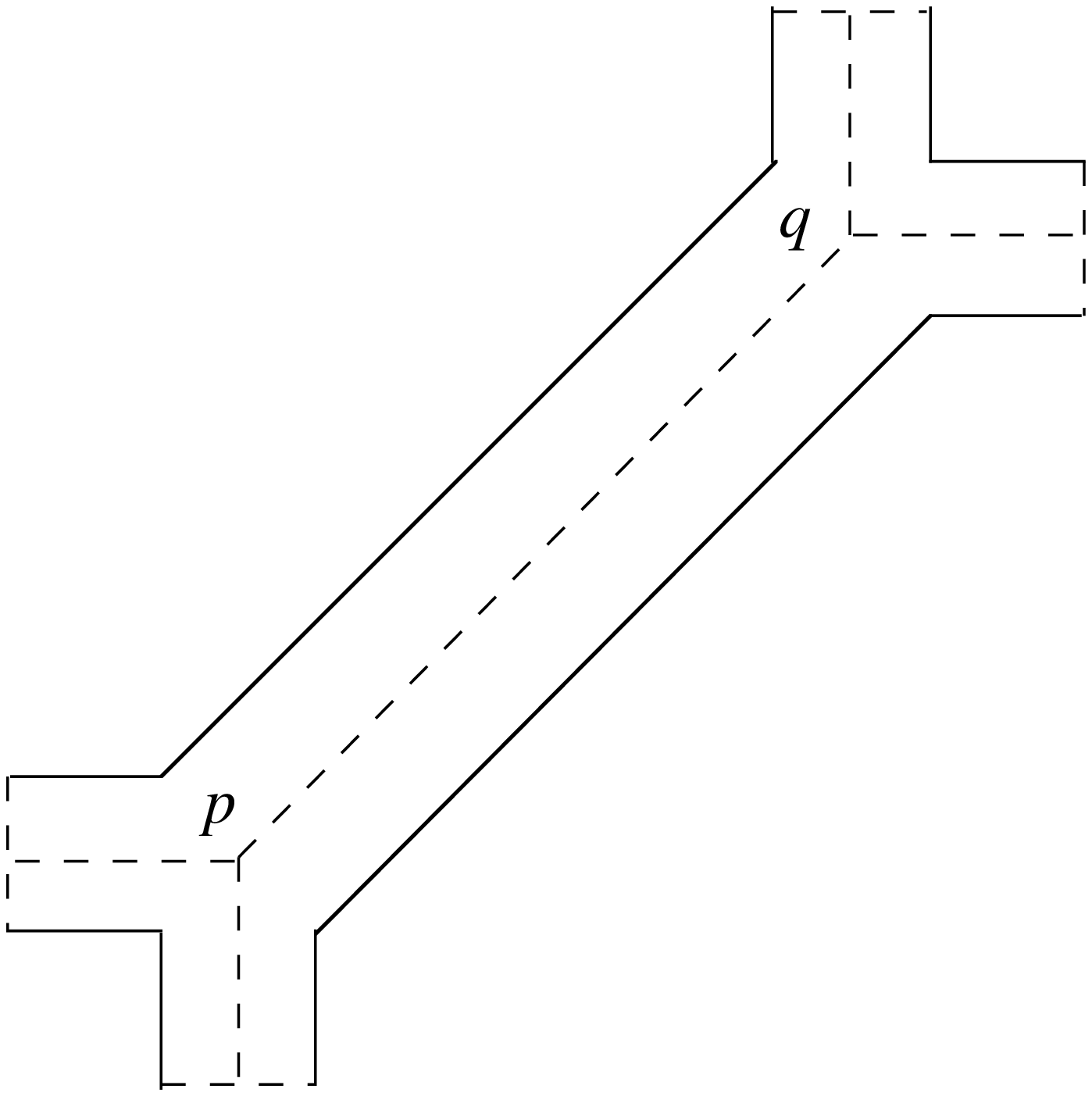}}
\vskip 20pt

Now consider the 2-manifold $M$ with boundary illustrated in the figure (with corners smoothed so that it has smooth boundary) --- in the figure, the top and bottom dashed lines, and the leftmost and rightmost dashed lines, are identified. It is clear that $M$ satisfies the geodesic condition of the Theorem, so the lower bound is valid for it. 
Assume that there is a vector field $A$ on $M$ having positive commutator with $H$. 
Notice that the two points $p$ and $q$ are such that there are three geodesics from $p$ to $q$: one in the direction $e_1 + e_2$, one in the direction $-e_1$ and one in the direction $-e_2$. Write $A = a_1 e_1 + a_2 e_2$.
Then \eqref{pos-def} implies that $a_1$ is increasing in direction $e_1$, and $a_2$ is increasing in direction $e_2$. Hence $a_1(p) > a_1(q)$, and $a_2(p) > a_2(q)$. On the other hand, $a_1 + a_2$ is increasing in direction $e_1 + e_2$, so this yields $a_1(p) + a_2(p) < a_1(q) + a_2(q)$, which is a contradiction. Hence no such vector field $A$ exists. 

\

\end{document}